\input amstex
\documentstyle{amsppt}

\centerline{\bf ON (SUB)STOCHASTIC AND TRANSIENT WEIGHTINGS}

\smallskip

\centerline{\bf OF INFINITE STRONG DIGRAPHS }

\smallskip

\centerline{\bf S.V. Savchenko}

\smallskip

\centerline{L.D. Landau Institute for Theoretical Physics, Russian Academy of Sciences}
\centerline{Kosygin str. 2, Moscow 119334, Russia}

\smallskip

\centerline{{\sl e-mail}:\ savch$\@$itp.ac.ru}

\bigskip

\centerline{\sl Dedicated to Richard Brualdi on the occasion of his 85th birthday}

\bigskip

{\bf Abstract.}
In the present paper,
for a given (possibly, infinite) strongly connected digraph
${\Cal D},$ we consider the class $\Cal{S}_{<}({\Cal D})$ of all truthly substochastic
weightings of ${\Cal D}$
(here, the word "truthly" means that there exists a vertex whose out-weight is strictly less than $1$).
For a finite subdigraph $\Cal{F}$ of $\Cal{D}$ weighted by
$S\in {\Cal S}_{<}({\Cal D}),$ let $\ell_{max}(\Cal{F})$ be the length of its
longest directed cycle and $\lambda_{S}(\Cal{F})$ be the Perron root (spectral radius)
of its weighted adjacency matrix.
We prove that the infimum of
$\ell_{max}(\Cal{F})\bigl(1-\lambda_{S}(\Cal{F})\bigr)$
taken over all $\Cal{F}$ is positive
for every $S\in \Cal{S}_{<}({\Cal D})$
if and only if
$\Cal{D}$ admits a finite cycle transversal.
The result obtained provides general theorems
on the set ${\Cal T}({\Cal D})$ of transient weightings of ${\Cal D}.$
In particular, we present a theorem of alternatives for
finite approximations to elements of
${\Cal T}({\Cal D})$ and simply reprove V. Cyr's criterion for
${\Cal T}({\Cal D})$ to be empty.


\smallskip


\bigskip

\centerline{\bf \S 1. Introduction}

\smallskip

Let ${\Cal D}$ be a (possibly, infinite) digraph with vertex-set $V({\Cal D})$
and arc-set $A({\Cal D}).$ We assume that
${\Cal D}$ can have possible loops, but it admits no
multiple arcs. By definition, ${\Cal D}$
is {\sl strongly connected} (or, simply,
{\sl strong}) if for any pair of its distinct vertices $v$ and $w,$
there is a path from $v$ to $w$ in ${\Cal D}.$
For each vertex $v$ of a strong digraph,
there exists at least one arc leaving $v.$
If a finite strong digraph has the property that
exactly one arc leaves each of its vertices, then
it is a {\sl directed cycle}. Its {\sl length} is
the number of its arcs.  Every strong digraph $\Cal{D}$
of order $|V({\Cal D})|$ at least two admits a cycle $\gamma$ of length $\ell(\gamma)$ at least two.
If $\ell(\gamma) < |A({\Cal D})|,$ then we say that $\gamma$
is a {\sl proper} cycle of $\Cal{D}$ and write $\gamma\subset \Cal{D}.$
The fact that $\gamma$ can coincide with $\Cal{D}$ is reflected as follows
$\gamma\subseteq \Cal{D}.$ The same rule also concerns arbitrary subdigraphs of $\Cal{D}.$

Throughout the paper, we consider
only directed cycles. In particular, we assume that
a loop is a directed cycle of length one and an undirected
edge is a directed cycle of length two. Denote by $\ell_{min}(\Cal{D})$
and $\ell_{max}(\Cal{D})$ the lengths of the shortest and longest cycles of $\Cal{D},$
respectively.  If for each integer $L>0,$ there exists a cycle of $\Cal{D}$ whose
length is not less than $L,$ then we assume that $\ell_{max}(\Cal{D})=\infty.$
A vertex-set is a {\sl cycle transversal} of ${\Cal D}$
if the vertex-deleted (sub)digraph
obtained by removing all its vertices from ${\Cal D}$
together with all their incident arcs
has no cycles, at all. Denote by $|sct(\Cal{D})|$ the size
of a smallest cycle transversal of $\Cal{D}.$ If $\Cal{D}$ admits no
finite cycle transversal, then we assume that $|sct(\Cal{D})|=\infty.$

A {\sl weighting} $W$ of ${\Cal D}$ is a positive function on its arcs.
The {\sl weighted adjacency matrix} $A_{W}(\Cal{D})$ is defined
as follows: its $(v,w)$th entry $A_{W}\bigl(\Cal{D}\bigr)(v,w)$
equals the weight of $(v,w)$ if $(v,w)$ is
an arc of the digraph, and coincides with $0,$ otherwise.
If $\Cal{D}$ is finite, then
{\sl the Perron root} $\lambda_{W}(\Cal{D})$ of $\Cal{D}$ weighted by $W$
is the Perron root (spectral radius)
of the corresponding weighted adjacency matrix $A_{W}(\Cal{D})$.
In the infinite case, the quantity $\lambda_{W}(\Cal{D})$ is defined as
the supremum of $\lambda_{W}(\Cal{F})$ taken over all finite (strong)
subdigraphs $\Cal{F}$ of $\Cal{D}.$

A weighting is called {\sl substochastic} if for
each vertex of ${\Cal D},$ its out-weight (i.e.
the sum of the weights of the arcs leaving it)
is not greater than $1.$ It will be {\sl truthly} substochastic if
there exists at least one vertex such that
its out-weight is strictly less than $1.$
Denote by ${\Cal S}({\Cal D})$ the class of all
substochastic weightings of ${\Cal D}$ and by ${\Cal S}_{<}({\Cal D})$
the class of all truthly substochastic weightings of ${\Cal D}.$
In the present paper, we study
the infimum of $\ell(\gamma)\bigl(1-\lambda_{S}(\gamma)\bigr)$
taken over all proper cycles $\gamma$ of $\Cal{D}$
and
the infimum of
$\ell_{max}(\Cal{F})\bigl(1-\lambda_{S}(\Cal{F})\bigr)$
taken over all finite subdigraphs $\Cal{F}$ of $\Cal{D}.$

It is shown that the first quantity
is positive
on $\Cal{S}({\Cal D})$ if and only if $\Cal{D}$
admits a finite cycle transversal (see Proposition 1).
The same also holds for the second quantity  defined on $\Cal{S}_{<}({\Cal D})$
(see Proposition 3).
In particular, this implies that the strict inequality $\lambda_{S}(\Cal{D})<1$
holds for each $S\in \Cal{S}_{<}(\Cal{D})$ if $\Cal{D}$
admits a finite cycle transversal and the lengths of all its cycles
are uniformly bounded above. Proposition 2 shows that these two conditions
are also necessary. So, we have a criterion for positiveness
of $1-\lambda_{S}(\Cal{D})$ on $\Cal{S}_{<}(\Cal{D})$
(see Proposition 4).
The proof means that for arbitrary ${\Cal D},$
the propositions will be true if one replaces $\Cal{S}_{<}({\Cal D})$ by the class
$\Cal{SS}_{<}({\Cal D})$ of all {\sl strictly} substochastic weightings of ${\Cal D}.$
In this case, for every vertex of ${\Cal D},$ its out-weight
is strictly less than $1.$

The class $\Cal{S}(\Cal{D})$ can be also considered as the family of all
substochastic matrices whose digraph is ${\Cal D}.$ Recall that the digraph of a matrix
is defined as follows: its vertex-set coincides with the index-set and $(v,w)$
is its arc if and only if the $(v,w)$th entry is non-zero. A matrix is called {\sl irreducible}
if its digraph is strongly connected. The results of the
last section will be presented in purely matrix terms whose definitions
are given below (see also [9]).

For an arbitrary infinite non-negative matrix $M,$ let $\lambda_{n}(M)$ be
the supremum of the spectral radii of its principal submatrices of order $n$
and $\lambda(M)$ be the limit of $\lambda_{n}(M)$ as $n\to\infty$ (i.e. $\lambda(M)$
is the supremum of the spectral radii of its finite principal submatrices).
The quantity $\lambda(M)$ is often called the {\sl intrinsic spectral
radius} of $M.$ We always assume that $\lambda(M)<\infty.$
By the Seneta-Sarymsakov theorem (see [11]),
the series $\sum\limits_{p=0}^{\infty}M^{p}z^{-p}$ is (entrywise) convergent
for each $z>\lambda(M).$ If this also holds for $z=\lambda(M),$ then
$M$ is {\sl transient}. Otherwise, it is {\sl recurrent}.

The notion of
recurrence is very important in probability theory of infinite Markov chains and
thermodynamic formalism of countable Markov shifts (see [4] and [7]). Here,
we mainly study the class ${\Cal T}({\Cal D})$
of all transient matrices whose digraph coincides with ${\Cal D}$
(or, all transient weightings of ${\Cal D}$) instead of
the class ${\Cal M}({\Cal D})$ of all
non-negative matrices $M$ with $\lambda(M)<\infty$
whose digraph is ${\Cal D}.$
However, the results obtained also imply sufficient conditions for
$M\in {\Cal M}({\Cal D})$ to be recurrent.

First of all, in Section 4, we show how V. Cyr's criterion
for ${\Cal T}({\Cal D})$ to be empty (see [3]) can be derived from
Proposition 4.
According to this criterion (in our formulation),
the set ${\Cal T}({\Cal D})$ is empty if and only if
$|sct(\Cal{D})|<\infty$ and $\ell_{max}(\Cal{D})<\infty$ (see our Theorem 1).
For this case, each $M\in {\Cal M}({\Cal D})$
has nice recurrent and analytic properties, while as it is shown in
Example 2, the rate of convergence of $\lambda_{n}(M)$ to $\lambda(M)$ can be
slower than any given negative power of $n.$

Note that finding relations between recurrent properties
of the original infinite non-negative matrix and the spectral
properties of its truncations is an old important problem (see [11]).
It is shown in Section 4 that Proposition 3
(see also Corollary 2) provides a
theorem (of alternatives)
for the quantity $n\bigl(\lambda(M)-\lambda_{n}(M)\bigr)$
on the class $\Cal{T}({\Cal D})$ (see Theorem 2).
The theorem states that {\bf either} its rate of convergence
to $0$ can be arbitrary fast {\bf or} it is always bounded
away from $0$ (the latter holds if and only if $|sct(\Cal{D})|<\infty$).
Obviously,  the same proposition will be true for the quantity
$n^{\alpha}\bigl(\lambda(M)-\lambda_{n}(M)\bigr)$ with any $\alpha \in [1,\infty).$
However, Example 1 presented in Section 2 implies that one cannot take
$\alpha$ in the open interval $(0,1),$ here.

\bigskip

\centerline{\bf \S 2. Positiveness of
$\inf\limits_{\gamma \subset \Cal{D}}
\ell(\gamma)\bigl(1-\lambda_{S}(\gamma)\bigr)$
on ${\Cal S}({\Cal D})$}

\smallskip

For a cycle $\gamma$ in ${\Cal D}$ and $S\in \Cal{S}(\Cal{D}),$
its {\sl weight} $S(\gamma)$ is the
product  of the weights of its arcs.
It is well known that the quantity $\lambda_{S}(\gamma)$
is equal to the $\ell(\gamma)$th root of $S(\gamma),$ i.e.
the {\sl gain} of $\gamma.$
So, the quantity
$\ell(\gamma)\bigl(1-\lambda_{S}(\gamma)\bigr)$
can be analyzed with the use of purely graph-theoretic methods only.

\smallskip

{\bf Proposition 1.} {\sl  For a strong digraph $\Cal{D},$
the infimum of $\ell(\gamma)\bigl(1-\lambda_{S}(\gamma)\bigr)$
taken over all proper cycles $\gamma$ of $\Cal{D}$ is positive on $\Cal{S}({\Cal D})$
if and only if
$\Cal{D}$ admits a finite cycle transversal:}
$$\inf\limits_{\gamma \subset \Cal{D}}
\ell(\gamma)\bigl(1-\lambda_{S}(\gamma)\bigr) > 0 \
\forall \  S\in \Cal{S}({\Cal D}) \iff \
|sct(\Cal{D})| < \infty.$$

{\bf Proof}. For a set (or sequence) $\Gamma$ of cycles, denote by
$V(\Gamma)$ and $A(\Gamma)$ the union of their vertices
and the union of their arcs, respectively.
Obviously, if $\{\gamma_{1},...,\gamma_{m}\}$
is a finite set of vertex-disjoint cycles in ${\Cal D},$
then either there exists a cycle $\gamma_{m+1}$
which has no common vertex with any of $\gamma_{k},$
where $k=1,...,m,$ or $\cup_{k=1}^{m}V(\gamma_{k})$ is a finite cycle
transversal of ${\Cal D}.$
Hence, if ${\Cal D}$
admits no finite cycle transversal, then there exists an infinite set
$\Gamma$ of vertex-disjoint cycles $\gamma_{k}$ of ${\Cal D}.$
For any sequence $\{c_{k}\}_{k=1}^{\infty}$ of positive numbers which are
strictly less than $1$,
assign the weight $c_{k}^{\frac{1}{\ell_{k}}},$
where $\ell_{k}$ is the length of $\gamma_{k},$ to each arc of $\gamma_{k}.$
Take any truthly substochastic extension of this weighting
from $A(\Gamma)$ to $A({\Cal D})$ and denote it by $S.$
Then $\lambda_{S}(\gamma_{k})=c_{k}^{\frac{1}{\ell_{k}}}.$ So, if
$\ell_{k}(1-c_{k}^{\frac{1}{\ell_{k}}})$ tends to $0$ as $k\to \infty,$
then the infimum of $\ell(\gamma)\bigl(1-\lambda_{S}(\gamma)\bigr)$
taken over all $\gamma \subset \Cal{D}$ equals $0.$

Assume now that ${\Cal D}$ admits a finite cycle transversal
$\{v_{1},...,v_{H}\},$ where $H<\infty.$
Take any stochastic $S\in \Cal{S}(\Cal{D}).$
Assume that for some $c\in (0,1),$
there are infinitely many cycles $\gamma$ of $\Cal{D}$ with
$S(\gamma)\ge c.$ As each $\gamma$ contains at least one vertex of
$\{v_{1},...,v_{H}\}$, for some $p\in \{1,...,H\},$
there are also infinitely many $\gamma$ with $S(\gamma)\ge c$
containing $v_{p}.$ Obviously, the sum of the weights of these $\gamma$ is infinite.
On the other hand, this sum is not greater than the probability
of returning to $v_{p}$ for the Markov chain whose transition matrix
is $A_{S}(\Cal{D})$ and hence, it does not exceed $1.$
This contradiction means that $0$ is a unique possible limit point of the weights of all cycles of $\Cal{D}.$
As $S(\gamma)<1$ for each $\gamma \subset \Cal{D}$,
the fact that $1$ is not a limit point of $S(\gamma)$ means that
the supremum of $S(\gamma)$ taken over all $\gamma \subset \Cal{D}$
is strictly less than $1.$
If this inequality holds for each stochastic weighting of ${\Cal D},$
then it also holds for each $S \in \Cal{S}({\Cal D}).$
From this fact and the evident inequality
$$1-S(\gamma)=\bigl(1-\lambda_{S}(\gamma)\bigr)\bigl(1+...+
\lambda_{S}(\gamma)^{\ell(\gamma)-1}\bigr)\le
\bigl(1-\lambda_{S}(\gamma)\bigr)\ell(\gamma)$$
it follows that the condition $|sct(\Cal{D})| < \infty$ is also sufficient
for positiveness of the infimum on $\Cal{S}({\Cal D}).$
The proposition is proved.

\smallskip

For $S\in {\Cal S}({\Cal D}),$
let $\omega_{S}(\Cal{D})$ be the supremum of $\lambda_{S}(\gamma)$
taken over all proper cycles $\gamma$ of ${\Cal D}.$
Proposition 1 allows us to formulate and prove a criterion for
positiveness of $1-\omega_{S}(\Cal{D})$
(i.e. the quantity $\inf\limits_{\gamma \subset \Cal{D}}
\bigl(1-\lambda_{S}(\gamma)\bigr)$ itself)
on $\Cal{S}({\Cal D}).$

\smallskip

{\bf Proposition 2.} {\sl  For a strong digraph $\Cal{D},$
the infimum of $1-\lambda_{S}(\gamma)$
taken over all proper cycles $\gamma$ of $\Cal{D}$ is positive on $\Cal{S}({\Cal D})$
if and only if
$\Cal{D}$ admits a finite cycle transversal and
the lengths of all cycles of ${\Cal D}$ are uniformly bounded above:}
$$\inf\limits_{\gamma \subset \Cal{D}}
\bigl(1-\lambda_{S}(\gamma)\bigr) > 0 \
\forall \  S\in \Cal{S}({\Cal D}) \iff \
|sct(\Cal{D})| < \infty
\and
\ell_{max}(\Cal{D})<\infty.
$$

{\bf Proof:}
Sufficiency of the conditions (i.e. $\omega_{S}(\Cal{D})<1$
for each $S\in {\Cal S}({\Cal D})$ if
$|sct(\Cal{D})| < \infty$ and
$\ell_{max}(\Cal{D})<\infty$) directly follows from Proposition 1.
So, it remains to show that
there exists $S\in {\Cal S}({\Cal D})$ with $\omega_{S}(\Cal{D})=1$
when $\ell_{max}(\Cal{D})=\infty.$
(As we have seen above, this always holds when $|sct(\Cal{D})|=\infty.$)

For an arbitrary decreasing
sequence of positive numbers $\epsilon_{k}<\frac{1}{2}$ convergent to zero,
let $\{\ell_{k}\}_{k=1}^{\infty}$ be a sequence of the lengths of
cycles in ${\Cal D}$ such that for each $k\ge 2,$
the length $\ell_{k}$ satisfies the inequalities
$\ell_{k}>L_{k-1}:=\ell_{1}+...+\ell_{k-1}$ and
$$\Bigl(\frac{1-\epsilon_{k}}{1-2\epsilon_{k}}\Bigr)^{\ell_{k}}>
\Bigl(\frac{2^{k}(1-\epsilon_{k})}{\epsilon_{k}}\Bigr)^{L_{k-1}}.\eqno(1)$$
It is possible under the condition $\ell_{max}(\Cal{D})=\infty$ because
$\ell_{k}$ can be arbitrary large and $\frac{1-\epsilon_{k}}{1-2\epsilon_{k}}>1$
(note that $\epsilon_{k}$ and $L_{k-1}$
are both already given).

Consider a sequence
$\Gamma=\{\gamma_{k}\}_{k=1}^{\infty}$ of (at the moment,
unweighted) cycles in ${\Cal D}$ such that the length of $\gamma_{k}$
is equal to $\ell_{k}.$
Let $\Gamma_{k-1}=\{\gamma_{1},...,\gamma_{k-1}\}.$
Assume that the weight of each arc in $\gamma_{1}$ equals
$1-\epsilon_{1}$ and then for each $k\ge 2,$
assign a weight equal to $1-\epsilon_{k}$ or $\frac{\epsilon_{k}}{2^{k}}$
to an arc $(v,w)$ of $\gamma_{k}$ when
$v\not\in V(\Gamma_{k-1})$ or $v\in V(\Gamma_{k-1}),$
but $(v,w)\not\in A(\Gamma_{k-1}),$ respectively.
If $(v,w)\in  A(\Gamma_{k-1}),$ then the weight of $(v,w)$ is already given.
Its value is equal to $1-\epsilon_{p}$ or $\frac{\epsilon_{p}}{2^{p}},$ where
$1\le p\le k-1.$ So, the weight of
each arc of $\gamma_{k}$ is not less than $\frac{\epsilon_{k}}{2^{k}}$
and at least $\ell_{k}-L_{k-1}$ arcs of $\gamma_{k}$ have weight equal to $1-\epsilon_{k}.$
By $(1)$ which can be rewritten as
$$\Bigl(\frac{\epsilon_{k}}{2^{k}}\Bigr)^{L_{k-1}}
\bigl(1-\epsilon_{k}\bigr)^{\ell_{k}-L_{k-1}} >
\bigl(1-2\epsilon_{k}\bigr)^{\ell_{k}},$$
the gain of $\gamma_{k}$ is not less than
$1-2\epsilon_{k}.$

Assume that
$v\in V(\Gamma)$ first appears in $\gamma_{k},$ i.e.
it is not contained in $\gamma_{1},...,\gamma_{k-1}$
but lies in $\gamma_{k}.$ Then the sum of the weights of all arcs
leaving $v$ and contained in $\Gamma$
does not exceed
$$1-\epsilon_{k}+\sum\limits_{p=k+1}^{\infty}\frac{\epsilon_{p}}{2^{p}}\le
1-\epsilon_{k}+\sum\limits_{p=k+1}^{\infty}\frac{\epsilon_{k}}{2^{p}}
=1-\epsilon_{k}+\frac{\epsilon_{k}}{2^{k}}\le 1-\frac{\epsilon_{k}}{2}<1.$$
So, the weighting constructed above is strictly substochastic
and any (sub)sto- chastic extension $S$
from $A(\Gamma)$ to $A({\Cal D})$ of it satisfies the condition $\omega_{S}(\Cal{D})=1.$
The proposition is proved.

\smallskip

For $n\ge \ell_{min}(\Cal{D})$ and $S\in \Cal{S}(\Cal{D}),$
we also consider the supremum $\omega_{S}(\Cal{D},n)$ of the gains
of proper cycles of length at most $n$ in $\Cal{D}$ weighted by $S.$
Obviously, $\omega_{S}(\Cal{D})$ is the limit of $\omega_{S}(\Cal{D},n)$
as $n\to \infty.$

\smallskip

{\bf Corollary 1.} {\sl
{\bf For} a given strong digraph ${\Cal D},$
{\bf either} for any positive function $g(n)>0,$ there exists
$S\in {\Cal S}_{<}({\Cal D})$ such that for
each $n\ge \ell_{min}(\Cal{D}),$
we have $1-\omega_{S}(\Cal{D},n)<g(n)$ {\bf or} for any
$S\in {\Cal S}({\Cal D}),$ there exists $c_{S}>0$ not
depending on $n$
such that $n\bigl(1-\omega_{S}(\Cal{D},n)\bigr)>c_{S}.$
The {\bf latter} alternative
holds if and only if ${\Cal D}$ admits a finite cycle transversal,
i.e. $|sct(\Cal{D})|<\infty.$
{\bf Finally}, for any $S\in {\Cal S}({\Cal D}),$
there exists $\delta_{S}>0$ not depending on $n$
such that $1-\omega_{S}(\Cal{D},n)>\delta_{S}$
if and only if $|sct(\Cal{D})|<\infty$ and $\ell_{max}(\Cal{D})<\infty.$}

\smallskip

{\bf Proof:} If $|sct(\Cal{D})|=\infty,$
then there is an infinite sequence $\Gamma=\{\gamma_{k}\}_{k=1}^{\infty}$
of vertex-disjoint cycles of $\Cal{D}.$
Denote by $\ell_{k}$ the length of the cycle $\gamma_{k}.$
The procedure of construction of $\Gamma$ presented in the beginning of
the proof of Proposition 1 implies that
we can chose $\Gamma$ so
that $\ell_{1}=\ell_{min}(\Cal{D})$ and
there exists  finite or infinite $p$ such that
$\ell_{k}<\ell_{k+1}$ for each $k\le p$ and
$\ell_{k}=L < \infty$ for any $k\ge p+1.$
\footnote[1]{In fact, we can assume that $p=0,1$ or $\infty.$
The latter means that the sequence $\{\ell_{k}\}_{k=1}^{\infty}$ is strictly increasing.}

Let us take a positive function $g(n)$ on ${\Bbb N}$ and then
construct $S\in {\Cal S}_{<}({\Cal D})$ such that
for each $n\ge \ell_{min}(\Cal{D})$, we have $1-\omega_{S}(\Cal{D},n)<g(n).$
Without loss of generality, we can assume that $g(n)<1$ for each $n\ge 1,$
the function $g(n)$ is strictly decreasing and tends to $0$ as
$n$ goes to infinity. Let us assign the weight
$1-g\bigl(\max\{\ell_{k+1},k\}\bigr)$ to each arc of $\gamma_{k}.$

Denote by $S$ some truthly substochastic extension to $A(\Cal{D})$
of the strictly substochastic weighting
of $A(\Gamma)$ presented above.
If $\ell_{k}\le n< \ell_{k+1},$ where $k\le p,$ then
$$1-\omega_{S}(\Cal{D},n)\le 1-\lambda_{S}(\gamma_{k})=1-\bigl(1-g(\ell_{k+1})\bigr)=
g(\ell_{k+1})< g(n).$$
In turn, for $n\ge \ell_{p+1}=L,$
we have $\omega_{S}(\Cal{D},n)=1-\lim\limits_{k\to\infty}g(k)=1$ and hence, for these values of $n,$
the inequality $0=1-\omega_{S}(\Cal{D},n)<g(n)$ holds.
The other statements of Corollary 1 directly follow from Propositions 1 and 2.
The corollary is proved.

\smallskip

{\bf Example 1.} Take an infinite path $1\to 2\to ...\to n \to ...$
and add arcs $(n,1),$ where $n=1,2,...,$ to it. Obviously, any cycle
of the digraph ${\Cal D}$ obtained in result has the form $1,2,...,n,1$
for some $n\ge 1.$ This means that the set $\{1\}$
is a smallest cycle transversal of ${\Cal D}$. To any $a\in (0,1)$ and each sequence of positive
numbers $f_{n}>0$ whose sum $f$ is not greater than $1,$ we can associate a
truthly substochastic
weighting $S$ of ${\Cal D}$ such that $S\bigl((1,1)\bigr)=af_{1},$
$S\bigl((1,2)\bigr)=1-f_{1},$ $S\bigl((n-1,n)\bigr)=
\frac{1-f_{1}-...-f_{n-1}}{1-f_{1}-...-f_{n-2}},$
and $S\bigl((n,1)\bigr)=\frac{f_{n}}{1-f_{1}-...-f_{n-1}}$ for $n\ge 2.$
In this case, we have
$S(\gamma_{n})=S\bigl((1,2)\bigr)S\bigl((2,3)\bigr)...
S\bigl((n-1,n)\bigr)S\bigl((n,1)\bigr)=f_{n}$ for the
cycle $\gamma_{n}$ with vertex-set $\{1,...,n\},$ where $n\ge 2.$

For $\epsilon>0,$ let us take
$$a_{\epsilon}=\sum\limits_{n=1}^{\infty}\frac{1}{n^{1+\epsilon}}\ \text{and}\
f_{n}=\frac{1}{a_{\epsilon}n^{1+\epsilon}}.$$
For $\epsilon^{\prime}>\epsilon,$
let $x_{n}=(1+\epsilon^{\prime}) \frac{\log n}{n}$ and
$\theta_{n}$ be the unique positive number such that
$(1+\theta_{n})^{n}=a_{\epsilon}n^{1+\epsilon}.$
Based on standard arguments of mathematical analysis, one can show
for sufficiently large $n,$ the inequality $(1+x_{n})^{n}>
a_{\epsilon}n^{1+\epsilon}$ holds and hence, $\theta_{n}\le x_{n}.$
Note that
$$1-\omega_{S}(\Cal{D},n)\le
1-\lambda_{S}(\gamma_{n})=1-\frac{1}{1+\theta_{n}}<\theta_{n}$$
and hence, $1-\omega_{S}(\Cal{D},n)=O(\frac{\log n}{n}).$
So, one cannot essentially improve the second
alternative in the statement of Corollary 1.

\bigskip

\centerline{\bf \S 3. Positiveness of
$\inf\limits_{\Cal{F}\subseteq\Cal{D}}
\ell_{max}(\Cal{F})\bigl(1-\lambda_{S}(\Cal{F})\bigr)$
on ${\Cal{S}_{<}({\Cal D})}$}

\smallskip

In this section, we prove similar results for the infimum of
$\ell_{max}(\Cal{F})\bigl(1-\lambda_{S}(\Cal{F})\bigr)$ taken over
all finite (not necessarily, induced) subdigraphs $\Cal{F}$ of $\Cal{D}.$
The necessity of their conditions directly follows from Propositions 1 and 2 which
will be also true if one replaces $\Cal{S}({\Cal D})$ by ${\Cal S}_{<}({\Cal D})$
in the statements. However, one cannot replace $\Cal{S}_{<}({\Cal D})$ by ${\Cal S}({\Cal D})$
in our propositions on the new quantity.

\smallskip

{\bf Proposition 3.} {\sl  For a strong digraph $\Cal{D},$
the infimum of $\ell_{max}(\Cal{F})\bigl(1-\lambda_{S}(\Cal{F})\bigr)$
taken over all finite (strong)
subdigraphs $\Cal{F}$ of $\Cal{D}$ is positive on $\Cal{S}_{<}({\Cal D})$
if and only if
$\Cal{D}$ admits a finite cycle transversal:}
$$
\inf\limits_{\Cal{F}\subseteq\Cal{D}}
\ell_{max}(\Cal{F})\bigl(1-\lambda_{S}(\Cal{F})\bigr) > 0\
\forall \  S\in \Cal{S}_{<}({\Cal D}) \iff \
|sct(\Cal{D})| < \infty.$$

{\bf Proof:}
Let $\Cal{F}$ be a finite subdigraph
of $\Cal{D}$ weighted by $S\in \Cal{S}_{<}({\Cal D})$
and $A_{S}(\Cal{F})$ be its weighted adjacency matrix.
By the Boyle-Handelman conjecture [1] proved by A. Goldberger
and M. Neumann in [6], the inequality
$\det\bigl(I-A_{S}(\Cal{F})\bigr)\le 1-\lambda_{S}(\Cal{F})^{r},$
where $r$ is the number of non-zero eigenvalues of $A_{S}(\Cal{F}),$ holds.
As $$1-\lambda_{S}(\Cal{F})^{r}=
\bigl(1-\lambda_{S}(\Cal{F})\bigr)\bigl(1+\lambda_{S}(\Cal{F})+...+\lambda_{S}(\Cal{F})^{r-1}\bigr)
\le \bigl(1-\lambda_{S}(\Cal{F})\bigr)r,$$
we also have
$$\det\bigl(I-A_{S}(\Cal{F})\bigr)\le r\bigl(1-\lambda_{S}(\Cal{F})\bigr).$$

The Coates determinant formula [2] which is a direct consequence
of the Leibniz formula for determinants (see [5]) and the fact that
any permutation is the product of independent cyclical permutations imply that
$$\det\bigl(I-zA_{S}(\Cal{F})\bigr)=1+\sum\limits_{U\in \Cal{U}(\Cal{F})}
(-1)^{n(U)}S(U)z^{\ell(U)},\eqno(2)
$$
where $\Cal{U}(\Cal{F})$ is the set of all unions $U$
of vertex disjoint cycles in $\Cal{F}$,
$n(U)$ is the number of cycles in $U,$
$S(U)$ is the product of their weights,
and $\ell(U)$ is the sum of their lengths.
As
$$\ell(U)\le \ell_{max}(\Cal{F})n(U)\le \ell_{max}(\Cal{F})|sct(\Cal{F})|,$$
the degree of $\det\bigl(I-zA_{S}(\Cal{F})\bigr)$ (and hence, the number $r$ of
non-zero eigenvalues of $A_{S}(\Cal{F})$) is not greater than the product of
$|sct(\Cal{F})|$ and $\ell_{max}(\Cal{F})$. This implies that
$$\det\bigl(I-A_{S}(\Cal{F})\bigr)\le |sct(\Cal{F})|\ell_{max}(\Cal{F})
\bigl(1-\lambda_{S}(\Cal{F})\bigr).$$

Assume now that $\Cal{D}$ admits a finite cycle transversal
$\{v_{1},...,v_{H}\}$ with $H<\infty.$ Without loss of generality, we can
assume that $\{v_{1},...,v_{h}\}=V(\Cal{F})\cap \{v_{1},...,v_{H}\}.$
For $k=1,...,h+1,$ let $\Cal{F}^{(k)}$ be the (sub)digraph obtained from $\Cal{F}$ by removing the vertices
$v_{1},...,v_{k-1}$ and all arcs incident to them
(here, we assume that $\Cal{F}^{(1)}=\Cal{F}$).
Then $A_{S}(\Cal{F}^{(k)})$ is the (sub)matrix obtained from $A_{S}(\Cal{F})$ by removing its
$v_{1},...,v_{k-1}$th columns and rows.
Since $\Cal{D}$ is strong
(so, $A_{S}(\Cal{D})$ is irreducible)
and $\lambda_{S}(\Cal{D})\le 1,$
we have $\lambda_{S}(\Cal{F}^{(k)})<1$ for each $k=1,...,h+1.$
Moreover, as $\Cal{F}^{(h+1)}$ is acyclic, $\lambda_{S}(\Cal{F}^{(h+1)})=0$.
In particular,  $\det\bigl(I-A_{S}(\Cal{F}^{(h+1)})\bigr)=1.$
So, Cramer's rule implies that
$$\frac{1}{\det\bigl(I-A_{S}(\Cal{F})\bigr)}=\prod\limits_{k=1}^{h}
\frac{\det\bigl(I-A_{S}(\Cal{F}^{(k+1)})\bigr)}
{\det\bigl(I-A_{S}(\Cal{F}^{(k)})\bigr)}=
\prod\limits_{k=1}^{h}
\bigl(I-A_{S}(\Cal{F}^{(k)})\bigr)^{-1}(v_{k},v_{k})=$$
$$\prod\limits_{k=1}^{h}
\sum\limits_{p=0}^{\infty}A_{S}\bigl(\Cal{F}^{(k)}\bigr)^{p}(v_{k},v_{k})\le
\prod\limits_{k=1}^{h}\sum\limits_{p=0}^{\infty}A_{S}(\Cal{F})^{p}(v_{k},v_{k})\le
\prod\limits_{k=1}^{h}\sum\limits_{p=0}^{\infty}A_{S}(\Cal{D})^{p}(v_{k},v_{k}).$$

It is well known
that
any irreducible substochastic matrix with at least one row sum strictly less
than $1$ is $1$-transient, i.e.
the series of its positive powers is (entrywise) convergent
(this follows directly from Theorem 1 [10]).
So, if we take the inverse value of the product
$\prod\limits\limits_{k=1}^{H}
\sum\limits_{p=0}^{\infty}A_{S}(\Cal{D})^{p}(v_{k},v_{k})<\infty$
as $d_{S},$ then we obtain $\det\bigl(I-A_{S}(\Cal{F})\bigr)\ge d_{S} >0$
and hence,
$$0< \frac{d_{S}}{H} \le
\frac{\det\bigl(I-A_{S}(\Cal{F})\bigr)}{|sct(\Cal{F})|}\le
\ell_{max}(\Cal{F})\bigl(1-\lambda_{S}(\Cal{F})\bigr).$$
The proposition is proved.

\smallskip

Propositions 2 and 3 taken together imply the following generalization
of the well-known fact that the spectral radius (Perron root) of a finite
irreducible truthly substochastic matrix is strictly less then one.
\footnote[2]{The proposition presented below can be also formulated as
follows: $\lambda_{S}(\Cal{D})<1$
for each $S\in {\Cal S}_{<}({\Cal D})$
$\ \iff \ |sct(\Cal{D})| < \infty \and \ell_{max}(\Cal{D})<\infty.$
One can show that under these conditions, the equality
$\lambda_{S}(\Cal{D})=1$ holds for each stochastic $S\in {\Cal S}({\Cal D}).$
This generalizes the trivial fact that the spectral radius (Perron root)
of a finite stochastic matrix is always equal to one.}

\smallskip

{\bf Proposition 4.} {\sl  For a strong digraph $\Cal{D},$
the infimum of $1-\lambda_{S}(\Cal{F})$
taken over all finite (strong)
subdigraphs $\Cal{F}$ of $\Cal{D}$ is positive on $\Cal{S}_{<}({\Cal D})$
if and only if
$\Cal{D}$ admits a finite cycle transversal and
the lengths of all cycles of ${\Cal D}$ are uniformly bounded above:}
$$
\inf\limits_{\Cal{F}\subseteq\Cal{D}}
\bigl(1-\lambda_{S}(\Cal{F})\bigr) > 0\
\forall \  S\in \Cal{S}_{<}({\Cal D}) \iff \
|sct(\Cal{D})| < \infty
\and
\ell_{max}(\Cal{D})<\infty.
$$

{\bf Remark 1}.
The last (elementary) part of the proof of Proposition 3
implies that if $H=|sct(\Cal{D})|<\infty,$ then for any $S\in \Cal{S}_{<}(\Cal{D}),$
there exists $d_{S}>0$ such that for each finite subdigraph
$\Cal{F}$ (of order $n$) of $\Cal{D}$ weighted by $S,$
the inequality $d_{S}\le \det\bigl(I-A_{S}(\Cal{F})\bigr)$ holds.
Moreover, if $L=\ell_{max}(\Cal{D})<\infty,$ at that,
then the degree of $det\bigl(I-zA_{S}(\Cal{F})\bigr)$ is not greater than $HL,$ i.e.
$A_{S}(\Cal{F})$ admits
at most $HL$ non-zero eigenvalues. As $|\lambda_{i}|<1$ for each
(non-zero) eigenvalue $\lambda_{i}$ of $A_{S}(\Cal{F}),$
we have
$$0< d_{S}\le |\det\bigl(I-A_{S}(\Cal{F})\bigr)|=
|1-\lambda_{1}|...|1-\lambda_{n}|\le
\bigl(1-\lambda_{S}(\Cal{F})\bigr)2^{HL-1}.$$
This lower bound is much worse than $0< d_{S}\le
\bigl(1-\lambda_{S}(\Cal{F})\bigr)HL$ obtained in the proof of Proposition 3.
However, it allows us to prove Proposition 4 without the use of
the Boyle-Handelman conjecture.

\smallskip

For $n\in \Bbb{N}$  and $S\in \Cal{S}(\Cal{D}),$
let $\lambda_{S}(\Cal{D},n)$ be the supremum of the Perron roots
of (induced) subdigraphs on (at most) $n$ vertices of $\Cal{D}$
weighted by $S.$
As $\omega_{S}(\Cal{D},n)\le \lambda_{S}(\Cal{D},n),$
Corollary 1 taken together with Propositions 3 and 4 implies the following proposition.

\smallskip

{\bf Corollary 2.} {\sl
{\bf For} a given strong digraph ${\Cal D},$
{\bf either} for any positive function $g(n)>0,$ there exists
$S\in \Cal{S}_{<}({\Cal D})$ such that for
each $n\ge \ell_{min}(\Cal{D}),$
we have $1-\lambda_{S}(\Cal{D},n)<g(n)$ {\bf or} for any
$S\in \Cal{S}_{<}({\Cal D}),$ there exists $c_{S}>0$ not
depending on $n$ such that $n\bigl(1-\lambda_{S}(\Cal{D},n)\bigr)>c_{S}.$
The {\bf latter} alternative
holds if and only if ${\Cal D}$ admits a finite cycle transversal.
{\bf Finally}, for any $S\in \Cal{S}_{<}({\Cal D}),$
there exists $\delta_{S}>0$ not depending on $n$
such that $1-\lambda_{S}(\Cal{D},n)>\delta_{S}$
if and only if ${\Cal D}$ admits a finite cycle transversal and
the lengths of all cycles of ${\Cal D}$ are uniformly bounded above.}

\smallskip

{\bf Remark 2}.
The second Keilson-Styan-Vermes theorem [8]
whose proof is much simpler than the one of the Boyle-Handelman conjecture
can guarantee that
$$\det\bigl(I-A_{S}(\Cal{F})\bigr)
\le 1-\lambda_{S}(\Cal{F})^{n}\le n\bigl(1-\lambda_{S}(\Cal{F})\bigr),$$
where $n$ is the order of $\Cal{F}$ (see also Appendix A.1).
As $0< d_{S}\le \det\bigl(I-A_{S}(\Cal{F})\bigr)$
when $S\in \Cal{S}_{<}(\Cal{D})$ and $|sct(\Cal{D})|<\infty,$
this implies the lower bound
$n\bigl(1-\lambda_{S}(\Cal{D},n)\bigr)>c_{S}$ with $c_{S}=d_{S}$
in Corollary 2. A (very short) proof of this bound with another (positive)
value of $c_{S}$ is actually given in Appendix A.1 (see Lemma A.2 therein).

\bigskip

\centerline{\bf \S 4. Finite approximation alternative theorem
for the class $\bold{{\Cal T}({\Cal D})}$}

\smallskip

In this section, based on Corollary 2 obtained
for the class ${\Cal S}_{<}({\Cal D}),$
we produce a similar result on the class ${\Cal T}({\Cal D})$
of transient matrices whose digraph is $\Cal{D}.$
Note that according to Theorem 1 [10], an irreducible
non-negative matrix $M$ is transient if and only if
there exists a positive vector $\vec{\xi}$
such that the entries of $M\vec{\xi}$ are not greater than the corresponding entries of $\lambda(M)\vec{\xi}$
and for at least one entry, the strict inequality holds.
Let $V_{\vec{\xi}}$ and $V_{\vec{\xi}}^{-1}$
be the diagonal matrices whose $(v,v)$th entries are equal to $\vec{\xi}(v)$
and $\vec{\xi}(v)^{-1},$ respectively.
Then $\lambda(M)^{-1}V_{\vec{\xi}}^{-1}MV_{\vec{\xi}}$
is a truthly substochastic matrix whose intrinsic spectral radius is equal to $1$
and whose digraph is the same as that of $M.$
This fact taken together
with the final proposition of the statement of Corollary 2
(see also Proposition 4) implies V. Cyr's
theorem (see Theorem 2.1 [3]) which was proved by him
in terms of thermodynamic formalism of countable
Markov shifts.

\smallskip

{\bf Theorem 1 [3].} {\sl
For a strong digraph ${\Cal D},$ the set ${\Cal T}({\Cal D})$
is empty if and only if
${\Cal D}$ admits a finite cycle transversal and the lengths of all cycles of ${\Cal D}$
are uniformly bounded above:}
$${\Cal T}({\Cal D})=\emptyset
\iff
|sct(\Cal{D})|<\infty \and \ell_{max}(\Cal{D})<\infty.$$

\smallskip

{\bf Remark 3}. The original condition of Theorem 2.1 [3]
(which is equivalent to that of our Theorem 1)
can be formulated in standard terms as follows: there exists
a finite vertex-subset such that the lengths of all {\bf walks} on
the corresponding vertex-deleted digraph are uniformly bounded above. For
a given strong digraph $\Cal{D}$ which does not satisfy this condition,
the construction of a transient matrix whose digraph is $\Cal{D}$ given in [3]
is very complicated and takes more than 20 pages. In our paper,
the existence of such a matrix almost directly follows from
the proofs of Proposition 1 (the case of $|sct(\Cal{D})|=\infty$)
and Proposition 2 (the case of $\ell_{max}(\Cal{D})=\infty$)
because the construction presented therein
yields a truthly substochastic weighting $S$ of $\Cal{D}$ with $\omega_{S}(\Cal{D})=1$
and hence, also a truthly substochastic
(transient) matrix $S$ with $\lambda(S)=1$ whose digraph is $\Cal{D}$.

\smallskip

In turn, the first part of the statement of Corollary 2 (see also Proposition 3)
implies a theorem of alternatives for finite approximations
to elements of $\Cal{T}({\Cal D})$ (its reformulation in terms of
thermodynamic formalism will be given in our paper under preparation
"Finite Markov subshifts of countable Markov shifts" which will be submitted
to {\sl Functional Analysis and its Applications}).

\smallskip

{\bf Theorem 2.} {\sl
{\bf For} a given strong digraph ${\Cal D},$
{\bf either} for any positive function $g(n)>0,$ there exists
$M\in {\Cal T}({\Cal D})$ such that for
any $n\ge 1,$
we have $\lambda(M)-\lambda_{n}(M)<g(n)$ {\bf or} for any
$M\in {\Cal T}({\Cal D}),$ there exists $c_{M}>0$ not
depending on $n$ such that
$n\bigl(\lambda(M)-\lambda_{n}(M)\bigr)>c_{M}.$
The {\bf latter} alternative
holds if and only if ${\Cal D}$ admits a finite cycle transversal.}
\footnote[3]{Here, we implicitly
assume that $\ell_{max}(\Cal{D})=\infty.$ In the opposite case,
by Theorem 1, we have ${\Cal T}({\Cal D})=\emptyset$ and hence,
nothing to talk about. Note also that if $|sct(\Cal{D})|=\infty,$
then by Corollary 2, for any $g(n)>0$ $\exists$
$S\in \Cal{S}_{<}({\Cal D})$ with $\lambda_{S}(\Cal{D})=1$ such that for
each $n\ge \ell_{min}(\Cal{D}),$ we have $1-\lambda_{S}(\Cal{D},n)<g(n).$
Take $c\in (0,1)$ so that $c<g(n)$ for each $n=1,...,\ell_{min}(\Cal{D}).$ Then $cS$
gives a transient matrix $M\in {\Cal T}({\Cal D})$ with
$\lambda(M)-\lambda_{n}(M)<g(n)$ for each $n\ge 1.$}

\smallskip

Theorem 2 shows that only in the case where $|sct(\Cal{D})|<\infty,$
we can expect that a fast convergence of $\lambda_{n}(M)$ to $\lambda(M)$
automatically implies nice recurrent properties of $M\in {\Cal M}({\Cal D}).$

\smallskip

{\bf Corollary 3.} {\sl Assume that the digraph of $M$ admits a finite cycle transversal
and $\lambda(M)-\lambda_{n}(M)=o(n^{-1}).$ Then $M$ is recurrent.}

\smallskip

According to Theorem 1,
if not only $|sct(\Cal{D})|<\infty,$ but also $\ell_{max}(\Cal{D})<\infty,$
then each $M\in {\Cal M}({\Cal D})$ is recurrent.
In our paper under preparation
"A note on analytic properties of first return
generating functions",
we show that in fact, all entries of
$\sum\limits_{p=0}^{\infty}M^{p}z^{p}$
are rational functions.
So,
analytic properties of $M$
are also very close to those of finite matrices.

After all that has been said about ${\Cal D}$
with $|sct(\Cal{D})|<\infty$ and $\ell_{max}(\Cal{D})<\infty,$
one might expect that
these conditions on ${\Cal D}$ imply a fast convergence of the difference
$\lambda(M)-\lambda_{n}(M)$ to zero for each $M\in {\Cal M}({\Cal D}).$
However, the following example shows that it is not so even in the simplest case.

{\bf Example 2}. To every infinite sequence of positive numbers $a_{k}$ with
$\sum\limits_{k=1}^{\infty}a_{k}^{2}<\infty,$ let us associate
a non-negative symmetric matrix $M$ with $I(M)={\Bbb N}$ and positive entries
$M(1,k+1)=M(k+1,1)=a_{k}$ for each $k\ge 1$.
Obviously, all directed cycles of the digraph of $M$
have length equal to $2$ and contain a common vertex $1.$
(Note that in the context of
undirected graph theory, on the contrary, the graph of $M$ has no cycles, at all and in fact,
it is an infinite star with center at $1.$)
It is also easy to see that if the numbers $a_{k}$ are ordered
in descending, then $M_{n}$ with the index-set $I(M_{n})=\{1,...,n\}$
has the largest spectral radius among all principal submatrices
of order $n$ in $M.$
By $(2),$ for $n\ge 2,$
we have $\det(I-zM_{n})=1-b_{n}^{2}z^{2},$ where $b_{n}^{2}=a_{1}^{2}+...+a_{n-1}^{2}.$
This means that $\lambda(M_{n})=b_{n}$ and hence, $\lambda(M)=b,$
where $b^{2}=\sum\limits_{k=1}^{\infty}a_{k}^{2}.$
Thus,
$$\lambda(M)-\lambda(M_{n})=b-b_{n}=\frac{b^{2}-b_{n}^{2}}{b+b_{n}}=
\frac{\sum\limits_{k=n}^{\infty}a_{k}^{2}}{b+b_{n}}.$$
It remains to notice that $b_{n}\to b$ as $n\to \infty$ and the tail
$\sum\limits_{k=n}^{\infty}a_{k}^{2}$ of the convergent series
$\sum\limits_{k=1}^{\infty}a_{k}^{2}$ can sufficiently slow converge
to zero. In particular, if $a_{k}=k^{-\frac{1+\epsilon}{2}}$ for some $\epsilon >0,$
then $\sum\limits_{k=n}^{\infty}a_{k}^{2}\asymp n^{-\epsilon}$ and hence,
$\lambda(M)-\lambda(M_{n})\asymp n^{-\epsilon}.$

\bigskip

\centerline{\bf Appendix A.1. Determinant inequalities for finite substochastic matrices}

\bigskip

In this short note, based on Cramer's rule and the similarity-invariance of the determinant $det$ and
the trace $tr,$ we present a lemma from which
the determinant inequality $n\bigl(1-\lambda(S_{n})\bigr)\ge \det(I-S_{n})$ directly follows
for the spectral radius $\lambda(S_{n})$ of a substochastic matrix $S_{n}$ of order $n.$
Note that for the case where $\lambda(S_{n})=1,$ the inequality is trivial because
by the Perron-Frobenius theorem,  $\lambda(S_{n})$ is an eigenvalue of $S_{n}$ and hence,
$\det(I-S_{n})=0.$ The case where $\lambda(S_{n})<1$ will be considered below.

\smallskip

{\bf Lemma A.1.} {\sl  Let $S_{n}$ be a substochastic matrix of order $n$
with $\lambda(S_{n}) <1.$ Then}
$$\frac{1}{1-\lambda(S_{n})}\le
tr(I-S_{n})^{-1}\le \frac{n}{\det(I-S_{n})}.$$

{\bf Proof:} Let
$\lambda_{1},...,\lambda_{n}$ be the eigenvalues of $S_{n}.$ Without loss of generality,
we can assume that $\lambda_{1}=\lambda(S_{n}).$ As $\Re(\lambda_{i})<1$ for each $i=1,...,n,$
we have
$$\frac{1}{1-\lambda(S_{n})}\le \frac{1}{1-\lambda(S_{n})}+
\sum\limits_{i=2}^{n}\frac{1-\Re(\lambda_{i})}{|1-\lambda_{i}|^{2}}
=\Re\Bigl(\sum\limits_{i=1}^{n}\frac{1}{1-\lambda_{i}}\Bigr)=$$
$$\Re\Bigl(tr(I-S_{n})^{-1}\Bigr)=tr(I-S_{n})^{-1}.$$

Denote by $I(S_{n})$ the index-set of $S_{n}.$
For each $v\in I(S_{n}),$ let $\bigl(I-S_{n}\bigr)^{-1}(v,v)$ be
the $(v,v)$th entry of $\bigl(I-S_{n}\bigr)^{-1}$ and $S_{n}^{(v)}$ be the (sub)matrix
obtained from $S_{n}$ by removing its $v$th column and row.
By Cramer's rule, we have
$$tr(I-S_{n})^{-1}=\sum\limits_{v\in I(S_{n})}\bigl(I-S_{n}\bigr)^{-1}(v,v)=
\sum\limits_{v\in I(S_{n})}
\frac{\det\bigl(I-S_{n}^{(v)}\bigr)}{\det\bigl(I-S_{n}\bigr)}
\le \frac{n}{\det(I-S_{n})}$$
because
$$\det(I-S_{n}^{(v)})=\exp\bigl(tr\ln(I-S_{n}^{(v)})\bigr)=
\exp\Bigl(-\sum\limits_{k=1}^{\infty}tr \frac{S_{n}^{(v)k}}{k}\Bigr)\le 1.$$
The lemma is proved.

\smallskip

The proof of Lemma A.1 also implies that
$$\frac{1}{n\bigl(1-\lambda(S_{n})\bigr)}\le
\max\limits_{v\in I(S_{n})}\bigl(I-S_{n}\bigr)^{-1}(v,v)\le \frac{1}{\det(I-S_{n})}.$$
In Section 3, based on Cramer's rule, we show that
$$\frac{1}{\det(I-S_{n})}\le
\prod\limits_{w\in W}\bigl(I-S_{n}\bigr)^{-1}(w,w),$$
where $W$ is a cycle transversal of the digraph $\Cal{D}_{n}$ of $S_{n}.$
(This allowed us to prove that if the digraph
of an infinite irreducible truthly substochastic matrix $S$ admits a finite cycle transversal,
then $n\bigl(1-\lambda_{n}(S)\bigr)$ is bounded away from zero.)
So, we have
$$\max\limits_{v\in I(S_{n})}
\bigl(I-S_{n}\bigr)^{-1}(v,v)\le
\prod\limits_{w\in W}\bigl(I-S_{n}\bigr)^{-1}(w,w).\eqno(A.1)$$
The following proposition yields a better upper bound on the maximum of
diagonal entries of $\bigl(I-S_{n}\bigr)^{-1}.$
Its proof is based on purely graph-theoretic arguments and does not use
Cramer's rule, at all.

\smallskip

{\bf Lemma A.2.} {\sl  If $W$ is a cycle transversal of the digraph $\Cal{D}_{n}$ of
a substochastic matrix $S_{n}$  (of order $n$) with $\lambda(S_{n})<1,$ then
for each $v\in I(S_{n}),$ we have}
$$\bigl(I-S_{n}\bigr)^{-1}(v,v)\le
1+\sum\limits_{w\in W}\sum\limits_{p=1}^{\infty}
S_{n}^{p}(w,w).$$

{\bf Proof:}
It is not difficult to check that
$$\bigl(I-S_{n}\bigr)^{-1}(v,v)=1+\sum\limits_{p=1}^{\infty}S^{p}_{n}(v,v)=
\sum\limits_{\gamma\in \Cal{CW}(\Cal{D}_{n},v)}S(\gamma),$$
where $\Cal{CW}({\Cal D}_{n},v)$ is the set of all closed walks
starting (and ending) at $v$ in the digraph $\Cal{D}_{n}$ of $S_{n}$
and $S(\gamma)$ is the weight of $\gamma,$
i.e. the product of the entries of $S_{n}$ along $\gamma$.
Obviously,
any closed walk can be contracted in
a cycle. So, any $\gamma\in \Cal{CW}({\Cal D}_{n},v)$ contains a vertex belonging to $W.$
For any such $\gamma=v_{0},...,v_{i},...,v_{\ell}$ (with $v_{0}=v_{\ell}=v$),
let $t$ be the smallest $i$ such that $v_{i}\in W$
(in particular, this implies that $v_{k}\neq v$ for $k=1,...,t-1$).
Define $f(\gamma)$ as a sequence
$v_{t},...,v_{\ell-1},v_{0},...,v_{t}\in \Cal{CW}({\Cal D}_{n},v_{t})$ of the same length $\ell$
(here, $v_{0}$ is the last member coinciding with $v$).
\footnote[4]{If we set $\gamma_{vw}=v_{0},...,v_{t}$ and
$\gamma_{wv}=v_{t},...,v_{\ell}$ (recall that $v_{0}=v_{\ell}=v$ and
$v_{t}=w$), then we can write
$f(\gamma_{vw}\gamma_{wv})=\gamma_{wv}\gamma_{vw}.$}
Then $f$ is a one-to-one map of $\Cal{CW}({\Cal D}_{n},v)$ into
$\bigcup_{w\in W}\Cal{CW}({\Cal D}_{n},w)$ (in other words,
$f(\gamma^{\prime})\neq f(\gamma)$ if $\gamma^{\prime}\neq \gamma$).
As $S\bigl(f(\gamma)\bigr)=S\bigl(\gamma\bigr),$
we have
$$\bigl(I-S_{n}\bigr)^{-1}(v,v)=1+\sum\limits_{\gamma\in \Cal{CW}(\Cal{D}_{n},v)}
S\bigl(f(\gamma)\bigr)\le$$
$$1+
\sum\limits_{w\in W}\sum\limits_{\gamma\in \Cal{CW}(\Cal{D}_{n},w)}
S(\gamma)=
1+\sum\limits_{w\in W}\sum\limits_{p=1}^{\infty}
S_{n}^{p}(w,w).$$
The lemma is proved.

\smallskip

{\bf Remark A.1.} The proof presented above is based on
the considering of closed walks on the digraph of $S_{n}.$
In fact, one can also give a proof of Cramer's identity
$$\bigl(I-S_{n}\bigr)^{-1}(v,v)=
\frac{\det\bigl(I-S_{n}^{(v)}\bigr)}{\det\bigl(I-S_{n}\bigr)}$$
which has the same nature. This directly follows from the fact that
$\bigl(I-zS_{n}\bigr)^{-1}(v,v)$ coincides with the zeta function associated with
all closed walks on $\Cal{D}_{n}$ containing the vertex $v,$
$\frac{1}{\det\bigl(I-zS_{n}^{(v)}\bigr)}$ is the zeta function associated with
all closed walk on $\Cal{D}_{n}-v,$ and $\frac{1}{\det\bigl(I-zS_{n}\bigr)}$
is the zeta function associated with all closed walks on $\Cal{D}_{n}$
(see Theorem 8.2 [7]).

\smallskip

In fact, Lemma A.2 implies that the product over elements of $W$ in the RHS of $(A.1)$
can be replaced by the elementary symmetric polynomial $\sigma_{k}$ of degree $k$
in $|W|$ variables for each $k=1,...,|W|.$ In other words,
for any cycle transversal $W=\{w_{1},...,w_{h}\}$
of the digraph $\Cal{D}_{n}$ of $S_{n}$ and each $k=1,...,h,$ we have
$$\bigl(I-S_{n}\bigr)^{-1}(v,v)\le
\sigma_{k}\Bigl(\bigl(I-S_{n}\bigr)^{-1}(w_{1},w_{1}),...,
\bigl(I-S_{n}\bigr)^{-1}(w_{h},w_{h})\Bigr).$$
This partially confirms the conjecture that any cycle transversal of $\Cal{D}_{n}$ always contains
a vertex at which the maximum of $\bigl(I-S_{n}\bigr)^{-1}(v,v)$ taken over all vertices
of $\Cal{D}_{n}$ is attained.

\bigskip

\centerline{\bf References}

\smallskip

[1]  M. Boyle and D. Handelman, The spectra of non-negative matrices
via symbolic dynamics, {\sl Annals of Math.} {\bf 133} (1991), 249-316.

\smallskip

[2] C. Coates,
Flow-graph solutions of linear algebraic equations,
{\sl IRE Transactions on circuit theory} {\bf 6} (1959), 170-187.

\smallskip

[3]  V. Cyr, Countable Markov shifts with transient potentials,
{\sl Proc. London Math. Soc.} {\bf 103} (2011), 923-949.

\smallskip

[4]  W. Feller, {\sl An introduction to probability theory and its
applications}, Vol. 1, Wiley, New-York, 1950.

\smallskip

[5]  F.R. Gantmacher, {\sl Matrix theory}, Chelsea, New York, 1959.

\smallskip

[6]  A. Goldberger and M. Neumann, An upper bound on the characteristic
polynomial of a nonnegative matrix leading to a proof of the Boyle-Handelman
conjecture, {\sl Proc. Amer. Math. Soc.} {\bf 137} (2009), 1529-1538.

\smallskip

[7] B.M. Gurevich and S.V. Savchenko, Thermodynamic formalism for
countable symbolic Markov chains, {\sl Uspekhi Mat. Nauk}  {\bf 53} (2)
(1998), 3-106; English transl. in  {\sl Russian Math. Surveys} {\bf 53} (2)
(1998), 245-344.

\smallskip

[8]  J. Keilson and G. Styan, Markov chains and M-matrices:
inequalities and equalities,
{\sl J. Math. Anal. Appl.} {\bf 41} (1973), 439-459.

\smallskip

[9] B. Mohar and W. Woess,
A survey on spectra of infinite graphs,
{\sl Bull. London Math. Soc.} {\bf 21} (1989), 209-234.

\smallskip

[10] W.E. Pruitt, Eigenvalues of non-negative matrices,
{\sl Ann. Math. Statist.} {\bf 35} (1964), 1797-1800.

\smallskip

[11] E. Seneta, Finite approximations to infinite non-negative matrices,
{\sl Math. Proc. Camb. Phil. Soc.} {\bf 63} (1967), 983-992.

\end{document}